\documentclass{amsart}
\usepackage{amsmath,amsthm}
\usepackage{amsfonts,amssymb}

\newtheorem{thm}{Theorem}[section]
\newtheorem{cor}[thm]{Corollary}
\newtheorem{lem}[thm]{Lemma}
\newtheorem{prop}[thm]{Proposition}

\theoremstyle{remark}

 \def\tb{{\mathbf t}}

 \def\xb{{\mathbf x}}

 \def\CR{{\mathcal R}}

 \def\CV{{\mathcal V}}

 \def\NN{{\mathbb N}}
 \def\PP{{\mathbb P}}
 
 \def\RR{{\mathbb R}}
 
        \def\sspan{\operatorname{span}}

\newcommand{\wt}{\widetilde}
\newcommand{\wh}{\widehat}

\begin{document}

\title
{Reconstruction of a polynomial from its Radon projections} 
 
\author{Borislav Bojanov and Yuan Xu}
\address{Department of Mathematics\\ University of Sofia\\
 Blvd. James Boucher 5\\1164 Sofia, Bulgaria}
\email{boris@fmi.uni-sofia.bg}

\address{Department of Mathematics\\ University of Oregon\\
    Eugene, Oregon 97403-1222.}
\email{yuan@math.uoregon.edu}

\date{\today} 
\thanks{The first author was supported by the Bulgarian Ministry of Science
and Education under Project MM--1402/2004. The second author was supported 
in part by the National Science Foundation under Grant DMS-0201669}

\begin{abstract}
A polynomial of degree $n$ in two variables is shown to be uniquely 
determined by its Radon projections taken over $[n/2]+1$ parallel lines 
in each of the $(2[(n+1)/2]+1)$ equidistant directions along the unit 
circle. 
\end{abstract}

\maketitle

\section{Introduction}
\setcounter{equation}{0}

Let $f$ be a function defined on the unit disk $B^2$ on the plane. A 
Radon projection of $f$ is the integral of $f$ over a line segment inside
$B^2$. More precisely, for any given pair $(\theta,t)$ of a real number
$t\in [-1,1]$ and any angle $\theta$, let  $I(\theta,t)$ denote the line 
segment inside the unit disk $B^2$, where the line passes through the
point $(t\cos \theta, t\sin\theta)$ and is perpendicular to the vector
$(\cos \theta, \sin \theta)$. Then
\begin{equation}\label{eq:1.1}
\CR_\theta(f;t) := \int_{I(\theta,t)} f(x,y) dx dy 
\end{equation} 
defines the Radon projection of $f$ on the line segment $I(\theta,t)$. 

The Radon transform $f \mapsto \{\CR_\theta(f;t)\}$ associates to $f$ a
family of univariate functions of $t$, parameterized by $\theta$. The
problem of reconstructing $f$ from full or partial knowledge of 
$\{\CR_\theta(f;t)\}$ has been studied by many authors. It plays an essential 
role in the computer tomography. It is well-known that the set of
Radon projections $\{ \CR_\theta(f;t): 0 \le t \le 1, 0 \le \theta \le 2 \pi\}$
determines $f$ completely. Furthermore, it is known that if $f$ has compact 
support in $B^2$, then $f$ is uniquely determined by any infinite set of 
Radon projections \cite{So}. In practice, however, the data set is usually
finite. Thus, the main problem is to obtain a good approximation to the 
function from a large collection of its Radon projections. See \cite{N} for 
the background and detail discussions. Using a finite data set, one may 
determine a polynomial whose Radon projections agree with the given data. 
Such an approach has been considered in \cite{HaSo,LS,Marr} in connection 
with computer tomography. 

The present paper concerns with the problem of whether a polynomial of 
degree $n$ in  two variables can be uniquely determined from a set of 
$(n+1)(n+2)/2$ distinct Radon projections. The solution depends on the 
arrangement of the lines on which the projections are taken.  One solution was 
given early in 
\cite{Marr}, in which the Radon projections are taken over all possible line 
segments $[\xb_i,\xb_j]$ joining $\xb_i$ and $\xb_j$, where $\xb_0, \xb_1, 
\ldots, \xb_{n+1}$ are $n+2$ equally spaced points on the boundary of the unit 
disk $B^2$. 
In \cite{H1} (see also \cite{CMS}), a more general result was proved, in which 
$\xb_0,\xb_1, \ldots, \xb_{n+1}$ are $n+2$ distinct points on the boundary 
of a convex set and the result was extended in \cite{H2} to higher dimensions. 
Recently this question was considered in \cite{BG}, in which the lines are 
taken in $(n+1)$ directions, with one line in the first direction, two lines 
in the second direction, etc. This set has the drawback of lacking 
symmetry.

The main result of the present paper shows that a set of Radon projections 
taken over $2 \lfloor  n/2 \rfloor +1$ parallel lines on each of the 
$2\lfloor (n+1)/2 \rfloor +1$ equidistant 
directions can determine a polynomial of degree $n$ uniquely. The set of the
corresponding line segments  possesses a rotation symmetry. If $n =2 m$ is a 
fixed positive integer, then our method requires the Radon data in $2m+1$ 
directions, which are taken to be in equally spaced angles along the unit 
circle, and in each direction we take $m+1$ parallel lines, associated to
$t_0, t_1,\ldots,t_m$. We prove that for almost
all choices of $\{t_k\}$ a polynomial $P_{n}$ of degree $n=2m$ in two variables
is uniquely determined by the set of its Radon projections over these line 
segments; more precisely, for any given $f$, $P_n$ is uniquely determined 
such that 
\begin{equation}\label{eq:1.2}
 \CR_{\phi_j} (P_{n}; t_k) = \CR_{\phi_j} (f; t_k),  
   \quad 0 \le j \le 2m, \quad  0\ \le k \le m. 
\end{equation} 
The similar construction works for $n = 2m-1$. Several examples of the
sets of points $\{t_k\}$ that define regular interpolation are given. 
Furthermore, the polynomial 
$P_n$ can be easily computed. Thus, the result appears to offer a 
simple way to recover a polynomial from its Radon data. 

The number of conditions \eqref{eq:1.2} is equal to $(n+1)(n+2)/2$, which is 
the dimension of the space of bivariate polynomials of total degree at 
most $n$.  Hence, \eqref{eq:1.2} can be interpreted as an interpolation problem
by polynomials based on line integrals. As in the case of interpolation
based on points, this bivariate problem is not always regular. 
The only general configurations in the literature that hold for every integer
$n$ are that of Hakopian \cite{H2} mentioned above and the non-symmetric
configuration given recently in \cite{BG}. The result of the present paper
provides another family of examples. The proof is based on an observation that
recovering polynomials from their Radon projections can be reduced,
using a formula in \cite{Marr}, to a family of univariate interpolation
problems that uses certain special classes of algebraic polynomials. In the 
case of \cite{BG}, the strategy of solving the interpolation problem resembles 
the approach that uses the Bezout theorem for pointwise interpolation. The 
approach in the present paper follows the strategy of \cite{BX1}, which uses 
equally spaced points on circles for pointwise interpolation, that 
allows one to step beyond the limitation of the Bezout Theorem. 

The paper is organized as follows. In the following section we prove the 
existence and the uniqueness of the reconstructing problem. In Section 3 we 
show how to construct the polynomial and give an outline of the algorithm.

\section{The existence and the uniqueness of the solution}
\setcounter{equation}{0}

\subsection{Preliminary}
Let $B^2= \{(x,y): x^2+y^2 \le 1\}$ denote the unit disk on the plane.
Let $\theta$ be the angle in the polar coordinates 
$$
x= r \cos \theta, \quad y = r \sin \theta, \qquad r \ge 0, \quad
   0 \le \theta \le 2 \pi.
$$ 
A line $\ell$ whose slope is $- \cot \theta$ is defined by the equation 
$$
  \ell(x,y) : = x \cos \theta + y \sin \theta - t = 0,
$$
where $t$ is a real number. Clearly the line $\ell$ passes through the point
$(t\cos\theta, t\sin \theta)$ and is perpendicular to the vector 
$\xi = (\cos \theta, \sin \theta)$. Since, for a fixed $t$, the line 
corresponding to $(\theta,t)$ coincides with the line corresponding to 
$(\pi + \theta, -t)$, we could assume $\theta \in [0, \pi)$. Alternatively, 
we could assume that $\theta \in [0, 2\pi)$ and $t \ge 0$.

We will also use $\ell$ to denote the set of points on the line $\ell$
and introduce  the notation
$$
I(\theta, t) = \ell \cap B^2, \qquad 0 \le \theta < \pi, \quad -1 \le t \le 1,
$$
to denote the line segment of $\ell$ inside $B^2$. The points on $I(\theta,t)$
can be represented as follows: 
$$
x = t \cos \theta - s \sin \theta, \quad  y = t \sin \theta + s \cos \theta,
$$
for $s \in [-\sqrt{1-t^2}, \sqrt{1-t^2}]$.

The Radon projection of a function $f$ in the direction $\theta$ with
a parameter $t \in [-1,1]$ is denoted by $\CR_\theta (f;t)$,   
\begin{align*}
\CR_\theta(f;t) := & \int_{I(\theta,t)} f(x,y) dx dy \\
    = & \int_{-\sqrt{1-t^2}}^{\sqrt{1-t^2}} f(t \cos \theta -s \sin\theta,
    t \sin \theta + s \cos\theta)ds.
\end{align*}
In the literature it is also called an $X$-ray. It is clear that 
$\CR_\theta(f;t)=\CR_{\pi+\theta}(f;-t)$ so that we can assume, for
the sake of definiteness, either $0 \le \theta < \pi$ or $0 \le t \le 1$. 

\subsection{Polynomial bases}
Let $\Pi_n^2$ denote the space of polynomials of total degree $n$ in 
two variables, which has dimension $\dim \Pi_n^2 = (n+1)(n+2)/2$. If 
$P \in \Pi_n^2$ then 
$$
 P(\xb) = \sum_{k=0}^n \sum_{j=0}^k c_{k,j} x^j y^{k-j}, \qquad \xb = (x,y).
$$ 
Let $U_k$ denote the Chebyshev polynomial of the second kind, 
$$
 U_k (x) = \frac{\sin(k+1)\theta}{\sin \theta}, \qquad   x = \cos \theta.
$$ 
For $\xi = (\cos \theta, \sin \theta)$ and $\xb = (x,y)$, the ridge 
polynomial $U_k(\theta; \cdot)$ is defined by
$$
U_k(\theta;\xb) := U_k(\langle \xb, \xi \rangle) = 
   U_k(x \cos \theta + y \sin \theta).
$$
Clearly $U_k$ is an element of $\Pi_k^2$ and it is constant on every line that 
is perpendicular to $\xi$. The Radon projection of this function in any
direction can be easily computed. This is a result due to Marr \cite{Marr} 
and it plays a central role in our discussion below. 

\begin{lem} \label{lem:Marr} 
For each $t \in (-1,1)$, $0 \le \theta, \phi \le 2\pi$, 
$$
 \CR_\phi(U_k(\theta;\cdot); t) = \frac{2}{k+1} \sqrt{1-t^2} U_k(t)
    U_k(\cos (\phi - \theta)).
$$
\end{lem} 

The following useful relation follows easily from Marr's formula (\cite{LS}), 
\begin{equation} \label{eq:3.3}
 \frac{1}{\pi} \int_{B^2} U_k(\theta;\xb) U_k(\phi; \xb) d \xb 
    = \frac{1}{k+1} U_k(\cos (\phi - \theta)).  
\end{equation}
Let $\CV_n$ denote the space of orthogonal polynomials of degree $n$
on $B^2$ with respect to the unit weight function; that is, $P \in \CV_n$
if $P$ is of degree $n$ and 
$$
\int_{B^2} P(\xb) Q(\xb) d\xb = 0, \qquad \hbox{for all $Q \in \Pi_{n-1}^2$}.  
$$
For special choices of $\theta$ and $\phi$, the equation \eqref{eq:3.3}
becomes an orthogonal relation. Indeed,  for $k \in \NN$, let 
$$
\xi_{j,k} = (\cos \theta_{j,k}, \sin \theta_{j,k}), \qquad 
   \theta_{j,k} :=\frac{j\pi}{k+1},  \quad 0 \le j \le k.  
$$
For fixed $k$, the points $\cos \theta_{j,k}$, $1 \le j \le k$, are zeros 
of $U_k$. The ridge polynomials $U_k(\theta;\xb)$ have the remarkable 
orthogonal property that $U_k(\theta; \cdot) \in \CV_k$. 
Since the dimension of $\CV_k$ is $k+1$, it follows from \eqref{eq:3.3} that
the set 
$$
  \PP_k: =  \left\{ U_k(\theta_{j,k}; \xb): 0 \le j \le k \right\}
$$
is a basis for $\CV_k$. In particular, this shows that the set 
$\{\PP_k: 0 \le k \le n \}$ is a basis for $\Pi_n^2$. 
Together with Lemma \ref{lem:Marr}, this proves the following result: 

\begin{lem} \label{lem:2.2} 
Every polynomial $P_n \in \Pi_n^2$ can be written uniquely as 
\begin{equation}\label{eq:Poly}
  P_n(\xb) = \sum_{k=0}^n \sum_{j=0}^k c_{j,k} U_k(\theta_{j,k}; \xb).
\end{equation}
Furthermore, for each $\phi$ and $t$, 
\begin{equation}\label{eq:RPoly}
\CR_\phi(P_n; t) = \sqrt{1-t^2} \sum_{k=0}^n \frac{2}{k+1} U_k(t) 
     \sum_{j=0}^k  c_{j,k} U_k(\cos (\phi - \theta_{j,k})). 
\end{equation}
\end{lem}

There are several other explicit orthogonal bases for $\CV_n$; see, 
for example, \cite{DX}. If $\{Q_j^k: 0 \le k \le n \}$ is an orthogonal
basis of $\CV_n$, then it was shown in \cite{X00} that 
\begin{equation}\label{eq:trans}
  U_k(\theta_{j,k}; \xb) = \frac{1}{k+1}\sum_{l=0}^k
      Q_l^k (x,y)  Q_l^k (\cos \theta_{j,k},\sin\theta_{j,k}). 
\end{equation}
One explicit basis of $\CV_n$, denoted by $Q_{j,i}(\xb)$, is given in terms 
of polar coordinates as follows (cf. \cite{DX}), 
\begin{align*}
Q_{j,1}(x,y) & = h_{n-2j,1}P_j^{(0, n-2j)}(2 r^2-1) r^{n-2j} \cos (n-2j)\theta,
   \quad 0 \le 2j \le n, \\
Q_{j,2}(x,y) & = h_{n-2j,2}P_j^{(0, n-2j)}(2 r^2-1) r^{n-2j} \sin (n-2j)\theta,
   \quad 0 \le 2j \le n-1, 
\end{align*}
where $P_j^{(\alpha,\beta)}$ denotes the usual Jacobi polynomial, 
$h_{0,1} = 1/(n+1))$ and $h_{j,1} = h_{j,2} = 1/(2n+2)$ for $1 \le 2j\le n-1$.
Using this basis and restricting $\xb$ to the boundary of $B^2$, the equation
 \eqref{eq:trans} becomes 
\begin{align} \label{eq:U}
\begin{split} 
& U_{2k}(\cos (\phi - \theta_{j,2k})) = 1 + 2 \sum_{l=1}^k
 \left(\cos 2l \theta_{j,2k} \cos 2 l \phi +
       \sin 2l \theta_{j,2k} \sin 2 l \phi \right), \\
& U_{2k-1}(\cos (\phi - \theta_{j,2k-1}))  =
2 \sum_{l=1}^k \left(\cos (2l-1) \theta_{j,2k-1} \cos (2 l-1) \phi 
  \right .\\
& \qquad \qquad \qquad \qquad \qquad \qquad \qquad 
 + \left .  \sin (2l-1) \theta_{j,2k-1} \sin (2l-1) \phi \right).  
\end{split}
\end{align}
The above equations also follow from the elementary trigonometric identities
$$
 1+ 2 \sum_{j=1}^k \cos 2 j\theta  = 
    \frac{\sin (2k+1)\theta}{\sin\theta} 
\quad\hbox{and}\quad
   2 \sum_{j=1}^k \cos (2j-1) \theta = \frac{\sin 2k\theta}{\sin\theta}
$$
with $\theta$ replaced by $\phi - \theta_{j,2k}$ and $\phi - \theta_{j,2k-1}$,
respectively. These relations allow us to rewrite \eqref{eq:RPoly} in a form 
that is easier to work with. In the following we use $\lfloor x \rfloor$
to denote the integer part of $x$. 

\begin{lem} \label{lem:2.3} 
Let $P_n$ be given as in \eqref{eq:Poly}. Then 
\begin{align} \label{eq:2.2}
& \frac{\CR_\phi(P_n; t)}{\sqrt{1-t^2}} = 
 \sum_{l=0}^{\lfloor n/2 \rfloor} U_{2 l}(t) \left[a_{0,2 l} + 
 2 \sum_{j=1}^l \left (a_{j,2l} \cos 2j \phi  
          + b_{j,2l} \sin 2j \phi \right ) \right] \\
& \quad  +  \sum_{l=1}^{\lfloor (n+1)/2 \rfloor} U_{2 l-1}(t) 
  \left[2 \sum_{j=1}^l 
    (a_{j,2l-1} \cos (2j-1) \phi + b_{j,2l-1} \sin (2j-1) \phi) \right],
 \notag
\end{align}
in which the coefficients $c_{j,k}$ in \eqref{eq:Poly} and $a_{j,k}$,
$b_{j,k}$ are related by
\begin{equation} \label{eq:2.3}
  c_{j,2l} = \frac{1}{2} a_{0,2l} + 
     \sum_{k=1}^l \left (a_{k,2l} \cos 2k \theta_{j,2l} 
          + b_{k,2l} \sin 2k \theta_{j,2l} \right ) 
\end{equation}
for $0 \le j \le 2l$, $0 \le l \le \lfloor n/2 \rfloor$, and 
\begin{equation} \label{eq:2.4}
  c_{j,2l-1} =  \sum_{k=1}^l \left (a_{k,2l-1} \cos (2k-1) \theta_{j,2l-1} 
          + b_{k,2l-1} \sin (2k-1) \theta_{j,2l-1} \right )
\end{equation}
for $0\le j \le 2l-1$, $1 \le l \le \lfloor (n+1)/2 \rfloor$.
\end{lem}

\begin{proof}
We only prove the case of $n =2m$. Comparing \eqref{eq:RPoly} and 
\eqref{eq:2.2} shows that 
\begin{align*}
&\frac{1}{2l+1} \sum_{j=0}^{2l} c_{j,2l} U_{2l} (\cos (\phi - \theta_{j,2l})) 
 = \frac{1}{2} a_{0,2l} + \sum_{j=1}^l (a_{j,2l} \cos 2j \phi 
   + b_{j,2l} \sin 2j \phi), \\
& \frac{1}{2 l}
\sum_{j=0}^{2l-1} c_{j,2l-1} U_{2l-1} (\cos (\phi - \theta_{j,2l-1})) \\
& \qquad\qquad\qquad
=\sum_{j=1}^l (a_{j,2l-1} \cos (2j-1) \phi + b_{j,2l-1} \sin (2j-1) \phi). 
\end{align*}
Using the equations in \eqref{eq:U} we can rewrite the left hand side
in the form of the right hand side, which leads to the equations:
$$
a_{j, 2l} = \frac{2}{2l+1} \sum_{k=0}^{2l} c_{k,2l}
      \cos 2j \theta_{k,2l}, \qquad    
b_{j, 2l} = \frac{2}{2l+1} \sum_{k=0}^{2l} c_{k,2l} \sin 2j \theta_{k,2l} 
$$
for $1 \le j \le l$ and $2 a_{0,2l}$ satisfies the above formula with 
$j =0$, and 
 \begin{align*}
& a_{j, 2l-1} = \frac{1}{l} \sum_{k=0}^{2l-1} 
    c_{k,2l-1} \cos (2j-1) \theta_{k,2l-1}, \\
& b_{j, 2l-1} = \frac{1}{l}  \sum_{k=0}^{2l-1} 
       c_{k,2l-1} \sin (2j-1) \theta_{k,2l-1} 
\end{align*}
for $1 \le j \le l$. For fixed $l$, the above linear relations can 
be written in the matrix form. After a proper normalization, the coefficient 
matrix turns out to be an orthogonal matrix. For example, for the even 
indices, we have
$$
 \frac{1}{\sqrt{2l+1}} \left[ \begin{matrix}
  1       & 1           & \cdots & 1 \\
 \sqrt{2}  &\sqrt{2} \cos \frac{2\pi}{2l+1} & \cdots 
                   & \sqrt{2} \cos \frac{2 l \pi}{2l+1} \\
  0  &\sqrt{2} \sin \frac{2\pi}{2l+1} & \cdots 
                   & \sqrt{2} \sin \frac{2 l \pi}{2l+1} \\
  \vdots & \vdots & \cdots & \vdots \\
  \sqrt{2}  &\sqrt{2} \cos 2 l \frac{2\pi}{2l+1} & \cdots 
                   & \sqrt{2} \cos 2l \frac{2 l \pi}{2l+1} \\
  0  &\sqrt{2} \sin 2l \frac{2\pi}{2l+1} & \cdots 
                   & \sqrt{2} \sin 2l \frac{2 l \pi}{2l+1} 
 \end{matrix}\right] \cdot  \left[\begin{matrix} 
     c_{0,2l}\\ c_{1,2l}\\ \vdots \\ c_{2l,2l}\end{matrix}\right] 
    = \frac{1}{\sqrt{2l+1}}  \left[\begin{matrix} 
     a_{0,2l}\\ \frac{a_{1,2l}}{\sqrt{2}}\\ \frac{b_{1,2l}}{\sqrt{2}} 
        \\ \vdots \\\frac{a_{l,2l}}{\sqrt{2}}\\ \frac{b_{l,2l}}{\sqrt{2}} 
\end{matrix}\right]. 
$$
The $(l+1)\times (l+1)$ matrix in the left hand side, including the 
factor $1/\sqrt{2l+1}$, is orthogonal. Hence, the linear system of equations 
can be easily reversed, which leads to the stated relations \eqref{eq:2.3} 
and \eqref{eq:2.4}. The formulas can also be verified directly by  
inserting \eqref{eq:2.3} and \eqref{eq:2.4} into the equations of 
$a_{j,k}$ and $b_{j,k}$ and using the 
well-known trigonometric identities such as 
$$
 \frac{1}{m+1} \sum_{k=0}^{m} \cos 2 j \theta_{k,m}  + i 
   \frac{1}{m+1} \sum_{k=0}^{m} \sin 2 j \theta_{k,m}  = 
    \frac{1}{m+1} \sum_{k=0}^{m} e^{2 i j \theta_{k,m}} =
  \delta_{0,j}
$$
for $0\le j \le m$, where $\delta_{0,j} = 1$ if $j \ne 0$ and 
$\delta_{0,j} = 0$ otherwise, and the elementary trigonometric identities
such as $2 \cos \theta \cos \phi = \cos (\theta + \phi) + \cos (\theta - \phi)$. 
\end{proof}

To reconstruct the polynomial, we will determine the coefficients 
$a_{j,k}$ and $b_{j,k}$ and use them in \eqref{eq:2.3} and \eqref{eq:2.4}
to determine $c_{j,k}$ in \eqref{eq:Poly}. The following representation
is useful. 

\begin{lem}  \label{lem:2.4} 
Let $n =2m$ or $n =2m-1$. Let $P_n$ be given as in \eqref{eq:Poly}. Then 
\begin{equation} \label{eq:AB}
 \frac{\CR_\phi(P_n; t)}{\sqrt{1-t^2}} =  A_0(t) + 
    \sum_{j=1}^n \left[ A_j(t) \cos j \phi +  B_j(t) \sin j \phi \right],
\end{equation}
where 
$$
 A_0(t) = \sum_{l=0}^{\lfloor n/2 \rfloor} a_{0,2 l} U_{2 l}(t), 
$$
and for $1 \le j \le m$, 
\begin{align} 
 A_{2j}(t) & =2 \sum_{l=j}^{\lfloor n/2 \rfloor} a_{j,2 l} U_{2 l}(t), \qquad  
 B_{2j}(t)  = 2 \sum_{l=j}^{\lfloor n/2 \rfloor} b_{j,2 l} U_{2 l}(t),
   \label{eq:Aeven}  \\
 A_{2j-1}(t) &  = 2 \sum_{l=j}^{\lfloor n/2 \rfloor} 
       a_{j, 2l-1} U_{2l-1}(t), \qquad  
 B_{2j-1}(t)  = 2 \sum_{l=j}^{\lfloor n/2 \rfloor} b_{j,2 l-1} U_{2 l-1}(t).
 \label{eq:Aodd} 
\end{align}
\end{lem}

\begin{proof}
This comes from changing the order of summations in \eqref{eq:2.2}. 
\end{proof}

\subsection{Existence and uniqueness of the solution}
We now state the Radon projections from which the polynomial $P_n$ 
can be uniquely recovered. These are given by 
\begin{equation} \label{eq:data}
 \CR_{\phi_{j,m}} (P_n; t_k),
        \qquad 0 \le k \le \lfloor n/2 \rfloor, \quad 
     0 \le j \le 2 m, \quad m =  \lfloor (n+1)/2 \rfloor, 
\end{equation}
which consists of total $(\lfloor n/2 \rfloor +1) ((2\lfloor (n+1)/2 \rfloor+1)
= (n+1)(n+2)/2$ many projections, the same as the dimension of $\Pi_n^2$. The 
angles are chosen to be equidistant, 
\begin{equation} \label{eq:theta}
\Theta_m : = \left \{ \phi_{j,m} =  \frac{2 j \pi}{2m+1}: 
   \quad  0 \le j \le 2m \right \}.  
\end{equation}

The reason that we choose the equidistant angles lies in the lemma below, 
which plays a key role in our study. Such a lemma appears first in 
\cite{BX1} and has been used in \cite{BX2} and \cite{X03}. 

\begin{lem} \label{lem:2.5}
For $\phi \in \Theta_m$ and $P_n\in \Pi_n^2$ with $n =2m$ or $n=2m-1$,  
\begin{align*}
&  \frac{\CR_\phi(P_n; t)}{\sqrt{1-t^2}} =  A_0(t) \\
& \qquad  + \sum_{j=1}^m
   \left[ \big( A_j(t) + A_{2m-j+1}(t)\big) \cos j \phi 
   +\big( B_j(t) - B_{2m-j+1}(t)\big) \sin j\phi \right], 
\end{align*}
where we assume that $A_{2m} = B_{2m} = 0$ if $n =2m-1$.
\end{lem}

\begin{proof}
Using the expression in the previous lemma, the proof follows from the 
fact that 
$$
\cos (2m-j+1)\phi = \cos j\phi \quad \hbox{and}\quad \sin (2m-j+1)\phi
  = - \sin j\phi
$$ 
for $\phi \in \Theta_m$. 
\end{proof}
For the expression in this lemma, the variable $t$ is in $(-1,1)$. We can
also state the result for $\phi \in \wt \Theta_m = \{(2j+1)\pi/(2m+1): 
0 \le j \le 2m\}$, for which the sign before $A_{2m-j+1}$ and $B_{2m-j+1}$ 
will be reversed. It is easy to see that $\wt \Theta_m = \Theta_m+ \pi$ 
modulus $2\pi$. The two expressions are consistent with the fact
that $\CR_{\theta} (f;t) = \CR_{\pi+\theta}(f;-t)$. 

We are now ready to prove our main result, which shows that the 
polynomial $P_n$ can be uniquely determined by the data \eqref{eq:data}. 
We need the following function classes: For $n = 2m$ define
$$
  X_j(t) = \{U_{2m}(t),U_{2m-2}(t), \ldots, U_{2j}(t), U_{2m-1}(t), 
        U_{2m-3}(t), \ldots, U_{2m-2j+1}(t) \} 
$$
for $1 \le j \le m-1$. We also consider $X_j$ as a column vector in 
$\RR^{m+1}$ and regard 
$$
  \Xi_j(\tb) = \left[ X_j(t_0), X_j(t_1), \ldots,X_j(t_m) \right]
$$
as an $(m+1) \times (m+1)$ matrix. 

A set of points $\tb=\{t_0, t_1,\ldots,t_m\}$ in $(-1,1)$ is said to be
asymmetric if $t$ and $-t$ do not both belong to $\tb$. 

\begin{thm} \label{thm:2.6}
Let $\tb:=\{t_0, t_1, \ldots, t_m\}$ be a set of asymmetric distinct points 
in $(-1,1)$ such that the matrices $\Xi_1(\tb), \Xi_2(\tb), \ldots, 
\Xi_{m-1}(\tb)$ are all nonsingular, then the polynomial $P \in \Pi_{2m}^2$ 
is uniquely determined by the set \eqref{eq:data} of its Radon projections.
\end{thm}

\begin{proof}
In order to prove that there is a unique solution, it is sufficient to 
show that if $\CR_{\phi_{j,m}} (P;t_k)= 0$ for $0 \le j \le 2m$ and 
$0 \le k \le m$, then $P(\xb) \equiv 0$.  

Setting its left hand side zero, the expression in Lemma \ref{lem:2.5}
shows that 
$$
 A_0(t_k)
 + \sum_{j=1}^m \left[ \big( A_j(t_k) + A_{2m-j+1}(t_k)\big) \cos j \phi 
   +\big( B_j(t_k) - B_{2m-j+1}(t_k)\big) \sin j\phi \right] 
  = 0
$$
for $0 \le k \le m$ and $\phi \in \Theta_m$. The left hand side is 
a trigonometric polynomial of degree $m$ and it vanishes on $2m+1$ 
points, hence, it follows from the uniqueness of the trigonometric
interpolation that the coefficients have to be zero; that is,
\begin{align} \label{eq:A(t)}
\begin{split}
A_0(t_k) =0, \qquad &  A_j(t_k) + A_{2m-j+1}(t_k) =0, \\
& B_j(t_k) - B_{2m-j+1}(t_k) =0, \qquad 1 \le j \le m
\end{split}
\end{align}
for $0 \le k \le m$. Since $A_j$ and $A_{2m-j+1}$ have different
parity, the coefficients in $A_j(t) + A_{2m-j+1}(t)$ will not combine
and there are exactly $m+1$ terms in this polynomial, written as a 
linear combination of the Chebyshev polynomials $\{U_k\}$. Let
$Y_0(t) = \{U_{2m}(t), U_{2m-2}(t), \ldots, U_2(t), U_0(t)\}$ and let  
\begin{align*}
 Y_{2j}(t) = \{U_{2m}(t), U_{2m-2}(t), \ldots, U_{2j}(t), U_{2m-1}(t), 
    U_{2m-3}(t), \ldots, U_{2m-2j+1}(t) \}, \\
 Y_{2j-1}(t) = \{U_{2m-1}(t), U_{2m-3}(t), \ldots, U_{2j+1}(t), U_{2m}(t), 
    U_{2m-2}(t), \ldots, U_{2m-2j}(t) \}. 
\end{align*}
Furthermore, let $\Xi_{Y_j}$ denote the matrix $ \Xi_{Y_j} = 
\left[Y_j(t_0), Y_j(t_1), \ldots, Y_j(t_m) \right]$. The coefficients of the 
linear systems of equations in \eqref{eq:A(t)} are $\Xi_{Y_0}$, $\Xi_{Y_{2j}}$
for $1 \le j \le m/2$, and $\Xi_{Y_{2j-1}}$ for $1 \le j \le (m+1)/2$.  
Hence, in order to prove that \eqref{eq:A(t)} implies $A_0(t) \equiv 0$, 
$A_j(t) = B_j(t) \equiv 0$ for $1 \le j \le m$, we need to show that 
these matrices are all invertible. However, it is easy to see that 
we have the relation $\Xi_{Y_{2(m-j)-1}} = \Xi_{Y_{2 j}}$, which shows, in 
particular, that $\Xi_{Y_m} = \Xi_{Y_{m-1}}$. Thus, we can deal with 
$\Xi_{Y_{2j}}$ for $0 \le j \le m-1$. Furthermore, $Y_0$ contains only
even polynomials $U_{2 l}$, $0 \le l \le m$. Using the notation
$s = t^2$, it follows that $U_{2l}(\sqrt{s})$ is a polynomial of degree $l$, 
so that the matrix $\Xi_{Y_0}$ is always invertible as $\tb$ is a set of 
asymmetric distinct points. Thus, since $\Xi_{Y_{2j}}= \Xi_j(\tb)$, our 
assumption on $\tb$ implies that all matrices are invertible.
\end{proof}

A similar theorem holds in the case of $n = 2m-1$. First we need to define
for $n =2m-1$, 
$$
  X_j(t) = \{U_{2m-2}(t), U_{2m-4}(t), \ldots, U_{2j}(t), U_{2m-1}(t), 
    U_{2m-3}(t), \ldots, U_{2m-2j+1}(t) \} 
$$
for $1 \le j \le m-1$ and 
$$
  X_m(t) = \{U_{2m-1}(t), U_{2m-3}(t), \ldots, U_3(t), U_1(t)\}. 
$$
We will keep the notion $\Xi_j(\tb)$ for the matrix built upon from the 
columns of $X_j$, which is an $m \times m$ matrix. 

\begin{thm} \label{thm:2.7}
Let $\tb=\{t_0, t_1, \ldots, t_{m-1}\}$ be a set of asymmetric distinct 
points in $(-1,0)\cup(0,1)$ such that the matrices $\Xi_1(\tb), \Xi_2(\tb),
\ldots, \Xi_{m}(\tb)$ are all nonsingular, then the polynomial $P \in 
\Pi_{2m-1}^2$ is uniquely determined by the set \eqref{eq:data} of its 
Radon projections. 
\end{thm}

That the points $t_0, t_1, \ldots, t_m$ are all in $(-1,0)\cup(0,1)$ means 
that none of the points can be zero. In fact, since $U_{2k-1}$ is an odd
polynomial, $\Xi_m(\tb)$ contains only odd polynomials, which vanish at 
the origin. Otherwise, the proof of this theorem is similar to that of 
the previous theorem. In the case of $n =2m$, the set 
$X_0=\{U_{2m}, U_{2m-2}, \ldots, U_0\}$ does not appear in the conditions
of the theorem since the interpolation at distinct points in $[0,1)$ by 
$X_0$ is always regular. Such a reduction of conditions does not appears 
to happen for $n =2m-1$. 

For each fixed $j$, the determinant of $\Xi_j(\tb)$ can be consider as a 
polynomial function in $\tb= (t_0, t_1, \ldots, t_m)$, so that 
$\det \Xi_j(\tb) =0$ defines a hypersurface of dimension $m$. Hence, 
the determinant is not zero for almost all choices of $\tb \in \RR^{m+1}$. 
Evidently, the same holds true for $m+1$ matrices. Hence, we
have the following corollary:

\begin{cor} \label{cor:2.7}
For $n = 2m$ or $2m-1$ and for almost all choices of distinct points 
$t_0, t_1, \ldots, t_m$ in $(-1,1)$ for $n =2m$ or in $(-1,0)\cup(0,1)$ for
$n =2m-1$, the polynomial $P_n \in \Pi_n^2$ is uniquely determined
by the set \eqref{eq:data} of its Radon projections. 
\end{cor}

\subsection{Choices of $t_k$}
We have shown that almost all choices of $\tb=\{t_0, t_1,\ldots, t_m\}$
will lead to the unique solution of recovering the polynomial. One may ask
the question if {\it any} choice of $\tb \subset (-1,1)^{m+1}$ will work. The 
answer, however, is negative as the following example shows.

\medskip\noindent
{\bf Example: $m=2$}. By Theorem \ref{thm:2.6} we only have to choose 
$t_0,t_1,t_2$ in $(0,1)$ such that the set $X_2 =\{U_2,U_3,U_4\}$
has nonsingular determinant $\Xi_2(\tb)$. Recall that 
$$
  U_2(t) = 4 t^2 -1, \quad  U_3(t) = 8 t^3 - 4t, \quad 
  U_4(t) = 16 t^4 - 12 t^2 +1. 
$$
We can compute the determinant $\Xi_2(\tb)$ explicitly and the result is  
\begin{align*}
 \det \Xi_2(\tb) = & \, 32 \prod_{1 \le i < j \le 3} (t_i - t_j) \\
   & \times \left [ 8 t_1^2 t_2^2 t_3^2 + 4 t_1 t_2 t_3(t_1 + t_2 +t_3)
     + (2 t_1^2 -1)(2 t_2^2 -1)(2 t_3^2 -1) \right ].
\end{align*}
Since $t_i \in (0,1)$, the first two terms in the square bracket are
positive, while the third term could be negative. This shows, in particular,
that $\det \Xi_2(\tb) \ne  0$ if one of the $t_i$ is $\sqrt{2} /2$ or if two 
of $t_1,t_2,t_3$ are less than $\sqrt{2}/2$ and the other one is greater
than $\sqrt{2}/2$. However, $\det \Xi_2(\tb)$ can be zero for some choices of 
$t_0, t_1, t_2$. 

\medskip

On the positive side, we give two results that provides sets of points
that will ensure the uniqueness of the reconstruction. The first result 
uses the following theorem due to Obrechkoff \cite{O}. 

\begin{lem} 
Let $d\mu$ be a nonnegative weight function on an interval and let 
$P_0, P_1, \allowbreak P_2, \ldots,$ be the orthogonal polynomials with 
respect to $d\mu$. Denote by $\alpha_n$ the largest zero of $P_n$. Then the 
number of zeros of the polynomial $a_0 P_0 + a_1 P_2 + \ldots + a_n P_n$, 
where $a_0, a_1, \ldots, a_n$ are real numbers, in the interval $(\alpha_n, 
+\infty)$ is at most the number of sign changes in the sequence of the
coefficients $a_0, a_1, \ldots, a_n$. 
\end{lem}

The Chebyshev polynomials $U_0, U_1, U_2, \ldots$ are orthogonal with 
respect to the weight function $\sqrt{1-x^2}$ on $[-1,1]$, so that the
above lemma implies the following:

\begin{prop}
Let $\tb = \{t_0, t_1,\ldots,t_m\}$ be a set of numbers that satisfies 
\begin{equation} \label{eq:2.15}
   \cos \frac{\pi}{2m+1} < t_0 < t_1 < \ldots < t_m < 1.
\end{equation}
Then the matrices $\Xi_1(\tb), \ldots, \Xi_{m-1}(\tb)$ are all nonsingular. 
\end{prop} 

\begin{proof}
If one of the matrices, say $\Xi_j(\tb)$, were singular, there would be
a nonzero polynomial $P \in \sspan X_j$ that vanishes on $\tb$. This means
that the number of zeros of $P$ in $(\alpha_n, \infty)$ would be $m+1$.
However, each set $\Xi_j$ has cardinality $m+1$ so that the number of sign 
changes in the sequence of the coefficients of $P$ is at most $m$. 
This is a contradiction to the conclusion of the lemma.
\end{proof}

The set of points given in this proposition ensures the uniqueness of 
determining a polynomial by the set of its Radon projections. The condition 
\eqref{eq:2.15} implies that all points are clustered toward the one end
of the interval $[-1,1]$. This appears to be neither practical nor a good 
choice for computation. In fact, the numerical test indicates that some of 
the matrices $\Xi_j$ tend to have very large condition numbers. 

Our second positive result is more interesting. Here the points $\{t_k\}$ 
are based on the zeros of the Chebyshev polynomials $U_{2m}$. It is 
well-known that these zeros are given by 
$$
   \eta_{j,2m}: = \cos \frac{j \pi}{2m+1}, \qquad j =1,2, \ldots,2m. 
$$
We state the result for $n = 2m$ first. 

\begin{thm} \label{thm:positive}
Let $t_0$ be any point in $(-1,1)$ such that $U_{2m}(t_0) \ne 0$. Let 
$t_j = \eta_{2j,2m}$ for $j = 1, 2, \ldots,m$. Then the matrices
$\Xi_1(\tb), \Xi_2(\tb), \ldots, \Xi_{m-1}(\tb)$ in Theorem 2.6 are all 
nonsingular. Consequently, the polynomial $P \in \Pi_{2m}^2$ is uniquely 
determined by the set \eqref{eq:data} of its Radon projection. 
\end{thm}

\begin{proof}
It is easy to see that the set $\tb$ is asymmetric. Hence, according to 
Theorem 2.6, we only need to show that $\Xi_j(\tb)$ is invertible for each 
$j =1,2,\ldots,m-1$. Using the explicit expression of $U_{2m}(t)$, 
it is easy to see that 
\begin{equation} \label{eq:Ueven}
    U_{2m-2j}(t_k) = U_{2m-2j}(\eta_{2k,2m}) = - U_{2j-1}(\eta_{2k,2m})
    = - U_{2j-1}(t_k).
\end{equation} 
Indeed, since $\sin(2m+1) \eta_{2k,2m} = 0$ and $\cos(2m+1) \eta_{2k,2m} = 1$,
it follows from the addition formula that $\sin (2m-2j+1)\eta_{2k,2m} = 
-\sin 2j \eta_{2k,2m}$, from which the equation follows. 

Since $U_{2m}(t_1) = \ldots = U_{2m}(t_m) =0$, we evidently have 
$$
   \det \Xi_j(\tb) = U_{2m}(t_0) \det \wt \Xi_j(\tb), 
$$
where $\wt \Xi_j(\tb)$ is a submatrix of $\Xi_j(\tb)$ with the column
$X_j(t_0)$ and the row containing $U_{2m}(t_j)$ of $\Xi_j(\tb)$ removed. 
Using the equation \eqref{eq:Ueven}, we can replace all rows of 
$\wt \Xi_j(\tb)$ that contain even Chebyshev polynomials by rows that 
contain odd Chebyshev polynomials. More precisely, we replace the rows
$(U_{2i}(t_1), \ldots, U_{2i} (t_m))$ for $i = j, j+1,\ldots, m-1$ by
rows $(U_{2m-2i-1}(t_1), \ldots, U_{2m-2i-1} (t_m))$, respectively. 
The new matrix has all elements given in terms of the Chebyshev polynomials
of the odd degree; hence,
$$
   \det \wh \Xi_j(\tb) = \varepsilon \det (U_{2i-1}(t_k))_{i=1, k=1}^{m, ~m}
$$ 
where $\varepsilon = \pm 1$. Assuming now that this determinant is zero. 
Then there exists a non-zero polynomial, $Q$, of the form 
$$
Q(x) = \sum _{i=1}^{m} b_i U_{2i-1}(x),
$$
which vanishes at $t_1, \ldots , t_m$. Since $Q$  is odd, it vanishes also
at $-t_1, \ldots , -t_m$. Observe that the sets $(-t_1, \ldots , -t_m)$ and
$(t_1, \ldots , t_m)$ do not overlap. We conclude that the polynomial $Q$ of 
degree $2m-1$ vanishes at $2m$ distinct points and, consequently, it vanishes 
identically, a contradiction. 
\end{proof} 

A similar theorem holds for $n =2m-1$, which we state below. 

\begin{thm}
Let $t_j = \eta_{2j,2m}$ for $j = 1, 2, \ldots,m$. Then the matrices
$\Xi_1(\tb), \Xi_2(\tb), \allowbreak \ldots, \Xi_m(\tb)$ in Theorem 2.7 
are all nonsingular. Consequently, the polynomial $P \in \Pi_{2m-1}^2$ is 
uniquely determined by the set \eqref{eq:data} of its Radon projection. 
\end{thm}

The proof is similar to that of Theorem \ref{thm:positive} and it uses the 
equation
$$
  U_{2m-2j}(\eta_{2k-1,2m}) = U_{2j-1}(\eta_{2k-1,2m}),
$$
which can be verified using elementary trigonometric identities.

We have conducted some numerical tests for identifying other sets of $t_k$ for 
which the matrices $\Xi_1(\tb), \ldots, \Xi_{m-1}(\tb)$ are nonsingular, so 
that the reconstruction of a polynomial from the set \eqref{eq:data} of 
its Radon projections is unique. For $n \le 20$, it turns out that both 
the equidistant points in $(0,1)$ and the Chebyshev points in $(0,1)$ 
work out. See, however, the discussion at the end of the next section.

\section{Reconstruction of polynomials}
\setcounter{equation}{0}

In this section we set up the algorithm that can be used to 
compute $P_n$ from the Radon projections. We only consider the case
$n =2m$. Let $\gamma_{j,k}$ be the data 
\begin{equation}\label{eq:A1}
\gamma_{j,k} = \CR_{\phi_j}(f; t_k) \Big / \sqrt{1-t_k^2}, \qquad 
   0 \le j \le 2m, \quad 0 \le k \le m,   
\end{equation}
where $\phi_{j,m}$ are given as in \eqref{eq:theta} and $t_k$ are 
distinct numbers in $[0,1)$, chosen in advance, such 
that the linear systems of equations \eqref{eq:A3} and \eqref{eq:A4}
below have unique solutions for all $j$. 

Recall that the Lagrange interpolation by trigonometric polynomials is 
used in the proof of Theorem \ref{thm:2.6}. The Lagrange interpolation 
based on the points in \eqref{eq:theta} is given explicitly by 
(see, for example, \cite{Z})
$$
  L_n f (\phi) = \sum_{j=0}^{2m} f(\phi_{j,m}) \ell_j(\phi), \qquad
  \ell_j(\phi) = \frac{\sin (m + \frac{1}{2}) (\phi - \phi_{j,m})}
      { (2m+1) \sin \frac{1}{2} (\phi - \phi_{j,m})},
$$  
where we assume that $f$ is the function being interpolated; that is,
$L_n f (\phi_{j,m}) = f (\phi_{j,m})$ for $0 \le j \le 2m$. Using the
well-known formula 
$$
  1 + 2 \cos \phi + \ldots + 2 \cos m \phi = 
    \frac{\sin (m + \frac{1}{2}) \phi }
         {\sin \frac{\phi}{2}},
$$
we can write $L_n f$ in the standard form of a trigonometric polynomial, 
$$
  L_n f (\phi) = m_0(f) + 2 \sum_{j = 1} \left( m_j^C (f) \cos j \phi + 
        m_j^S (f) \sin j \phi \right)
$$
where 
\begin{align*}
& m_j^C(f) =\frac{1}{2m+1} 
 \sum_{l =0}^{2m} f(\phi_{l,m}) \cos j \phi_{l,m}\\
& m_j^S(f) =\frac{1}{2m+1}  \sum_{l =0}^{2m} f(\phi_{l,m}) \sin j \phi_{l,m}.
\end{align*}
For each $k$, $0 \le k \le m$, we will use the above formulas with 
$f(\phi_{j,m}) = \gamma_{j,k}$ and we shall write $m_{j,k}^C= m_j^C(f)$
and $m_{j,k}^S= m_j^S(f)$ for this particular $f$. This will allow us
to determine the values $A_0(t_k)$, $A_j(t_k) + A_{2m-j+1}(t_k)$ and 
$B_j(t_k) - B_{2m-j+1}(t_k)$ for $1 \le j \le m$. 

The next step is to fix $j$ and use the values of $A_0(t_k)$, or
$A_j(t_k) + A_{2m-j+1}(t_k)$, or $B_j(t_k) - B_{2m-j+1}(t_k)$ for 
$k = 0, 1, \ldots, m$ to determine the coefficients of $A_j$ and $B_j$. 
For this purpose we consider the following linear systems of equations: 
For $1 \le j \le m/2$,
\begin{equation}\label{eq:A3}
 \sum_{l=j}^m d_{j,2l} U_{2l} (t_k) +
     \sum_{l=m-j+1}^m d_{m-j+1,2l} U_{2l-1} (t_k) = m_{2j,k}, 
   \quad 0 \le k \le m,     
\end{equation}
and for $1 \le j \le (m+1)/2$,
\begin{equation}\label{eq:A4}
 \sum_{l=j}^m d_{j,2l-1} U_{2l-1} (t_k) +
     \sum_{l=m-j+1}^m d_{m-j+1,2l} U_{2l} (t_k) = m_{2j-1,k}, 
   \quad 0 \le k \le m.     
\end{equation} 
The quantities $m_{j,k}$ will be specified later. 
It is easy to see that the coefficient matrix of these systems of 
equations are $\Xi_j$ in Theorem \ref{thm:2.6}. There are total $m$ 
system of linear equations, each of size $(m+1) \times (m+1)$. Solving 
these equations with proper $m_{j,k}$ will determine the coefficients of 
$A_j$ and $B_j$. 

The last step is to use \eqref{eq:2.3} and \eqref{eq:2.4} to get $c_{j,k}$,
which are the coefficients of $P_n$ in \eqref{eq:Poly}. 

\medskip \noindent
{\bf Algorithm:}
\medskip

\par\noindent
{\bf Step 1}. 
For $1 \le l \le m$ and $0 \le k \le m$, compute
\begin{equation}\label{eq:A2}
 m_{l,k}^C = \frac{1}{2m+1}\sum_{j=0}^{2m}\gamma_{j,k} \cos l \theta_j  
  \qquad \hbox{and} \qquad
 m_{l,k}^S =\frac{1}{2m+1} \sum_{j=0}^{2m}\gamma_{j,k} \sin l \theta_j. 
\end{equation}

\noindent 
{\bf Step 2}. Solve the systems of equations \eqref{eq:A3} and \eqref{eq:A4} 
to get the coefficients $a_{j,k}$ and $b_{k,j}$ in \eqref{eq:2.2}: 

{\it Case 1}. 
Solve the $m$ systems of \eqref{eq:A3} and \eqref{eq:A4} for $m_{j,k} = 
m_{j,k}^C$ to get
$$
  a_{j,k} := d_{j,k}, \qquad 1 \le j \le m \quad 0 \le k \le 2 m.
$$

{\it Case 2}. 
Solve the $m$ systems of \eqref{eq:A3} and \eqref{eq:A4} for $m_{j,k} = 
m_{j,k}^S$ to get
\begin{align*}
&  b_{j,2 l} := d_{j,2l},  \quad 1 \le j \le m/2, \quad\hbox{and}\quad 
   b_{j,2l} := - d_{j,2l}, \quad \frac{m+2}{2} \le j \le m,  \\ 
&  b_{j,2l-1} :=  d_{j,2l-1}, \quad 1 \le j \le \frac{m+1}{2}, 
  \quad\hbox{and}\quad 
   b_{j,2 l-1} := - d_{j,2l-1}, \quad \frac{m+1}{2}\le j \le m,
\end{align*}
where $1 \le l \le m$. 

\noindent 
{\bf Step 3}. Substitute the output $a_{j,k}$ and $b_{j,k}$ in Step 2 
into \eqref{eq:2.3} and \eqref{eq:2.4} to get $c_{j,k}$. Then 
the polynomial $P$ is given by \eqref{eq:Poly}.

\medskip

The output of this algorithm is the polynomial $P$ that satisfies 
\eqref{eq:1.2}. The main computation appears to be the Step 2. Initial
numerical experiments indicate that the linear systems of equations
\eqref{eq:A3} and \eqref{eq:A4} are ill-conditioned. However, among the
set of points that we tested, the largest condition number of the matrices 
corresponding to the equidistant points $t_k = (k+1)/(m+2)$, $0 \le k \le m$,
of $(0,1)$ is smaller than that corresponding to the Chebyshev points $t_k = 
\cos (k+1) \pi/(2m+4)$, $0 \le k \le m$ or the points in Theorem 
\ref{thm:positive}.

As pointed out by one referee, it is perhaps not surprising that the 
matrices of the linear systems \eqref{eq:A3} and \eqref{eq:A4} are 
ill-conditioned. Solving these linear systems means inverting 
the Radon transform, while it is known that the Radon reconstruction 
problem is ill posed. Thus, for $m$ large, the algorithm may not be 
useful for practical computation. On the other hand, our result 
shows that it is possible to determine a polynomial of degree $n$ 
uniquely from a set of $(n+1)(n+2)/2$ Radon data that concurs with 
the parallel geometry. This appears to be of independent interest.

\end{document}